\documentclass{elsarticle}
\usepackage{latexsym}
\usepackage{amsfonts}
\usepackage{amsmath}
\usepackage{amssymb}
\usepackage{amsthm}
\usepackage{xypic}
\usepackage{mathrsfs}
\usepackage{etoolbox}
\usepackage[hang]{footmisc}
\usepackage{abstract}

\input xy
\xyoption{all}

\newtheorem{prop}{Proposition}
\newtheorem{lem}{Lemma}
\newtheorem{thm}{Theorem}

\newtheorem*{pf}{Proof}

\begin{document}

\renewcommand{\thefootnote}{\fnsymbol{footnote}}
\title{The Witt Ring of a Curve with Good Reduction over a Non-dyadic Local Field}
\author{Jeanne M. Funk\thanks{CORRESPONDING AUTHOR\newline Department of Mathematics, Engineering and Computer Science\newline LaGuardia Community College \newline 31-10 Thomson Ave., Room E218-E \newline Long Island City, NY 11101 \newline jfunk@lagcc.cuny.edu\newline  1-718-482-5698}}
\author{Raymond T. Hoobler\thanks{Department of Mathematics\newline The City College of New York\newline NAC 6/203A\newline Convent Ave at 138th Street\newline New York, NY 10031\newline rhoobler@ccny.cuny.edu}}


\maketitle

\begin{abstract}

In this work, we present a generalization to varieties and sheaves of the fundamental ideal of the Witt ring of a field by defining a sheaf of fundamental ideals $\tilde{I}$ and a sheaf of Witt rings $\tilde{W}$ in the obvious way. The Milnor conjecture then relates the associated graded of $\tilde{W}$ to Milnor K-theory and so allows the classical invariants of a bilinear space over a field to be extended to our setting using \'etale cohomology.  As an application of these results, we calculate the Witt ring of a smooth curve with good reduction over a non-dyadic local field.
\end{abstract}

\saythanks

\clearpage
\renewcommand{\thefootnote}{\arabic{footnote}}\setcounter{footnote}{0}
\section{Introduction}

The goals of this work are twofold.  First, we present a generalization to the category of varieties and sheaves of the fundamental ideal of the Witt ring and of some techniques associated with this ideal which have proven useful in describing the Witt ring  of a field.  Second, we calculate the Witt ring of a smooth projective curve with good reduction over a non-dyadic local field as an application of these techniques.

There has been a great deal of work done toward generalizing the fundamental ideal to rings, schemes, and varieties. 
~\cite{auel} is a good summary of results.  Our approach is to define a sheaf of fundamental ideals $\tilde{I}$ and a sheaf of Witt rings $\tilde{W}$.  For low dimensional cases, the global sections $\mathbb{I}^n:=\Gamma(X,\tilde{I}^n)$ form a filtration of the Witt ring $W(X)$ and we can study the Witt ring via the quotients $\mathbb{I}^n/\mathbb{I}^{n+1}$.  We will show that $\tilde{I}^n/\tilde{I}^{n+1}$ is isomorphic to the Zariski sheafifications of \'{e}tale cohomology with $\mu_2$ coefficients (this follows from the Milnor conjecture and most of the work has already been done by Kerz, Milnor, Orlov, Vishik, and Voevodsky, among others).  We then show that $\tilde{W}/\tilde{I}$, $\tilde{I}/\tilde{I}^2$ and $\tilde{I}^2/\tilde{I}^3$ are isomorphic to the Zariski sheaves $\mathbb{Z}/2\mathbb{Z}$, $\mathcal{O}^{\times}/\mathcal{O}^{\times 2}$, and ${ }_2Br$ respectively and that the classes of elements in $W(X)/\mathbb{I}$, $\mathbb{I}/\mathbb{I}^2$, and $\mathbb{I}^2/\mathbb{I}^3$ are determined by rank, signed discriminant, and Witt invariant respectively.

For a smooth geometrically connected projective curve $C$ with good reduction $C_k$ over a non-dyadic local field with uniformizing parameter $\pi$ we can show that $\mathbb{I}^n/\mathbb{I}^{n+1}$ is trivial when $n>2$.  Thus, using the generalizations outlined above, we can describe each element of the Witt ring in terms of its rank, signed discriminant, and Witt invariant.  We may equate each space of odd rank to $<-1>$ in $W(X)/\mathbb{I}$ and each even rank space to $<1,\delta>$ in $\mathbb{I}/\mathbb{I}^2$ where $\delta$ is the signed discriminant.  It remains to describe the contribution from $\mathbb{I}^2/\mathbb{I}^3$, which is determined by the class of the Clifford algebra in the order $2$ part of the Brauer group.  We show that, over the function field of $C$, ${ }_2Br(C)$ consists of distinct "quaternions" $\big(\frac{\xi,\pi}{k(C)}\big)$, allowing us to equate each member of $\mathbb{I}^2$ to the four dimensional form $<1,-\xi,-\pi,\pi\xi>$ in $\mathbb{I}^2/\mathbb{I}^3$.  In this way, we can explicitly represent each element of the Witt ring of a sum of line bundles and prove that $W(C)\simeq W(C_k)[\mu_2]$, which nicely generalizes the classical result for the Witt ring of a local field ~\cite{lam}.

We would like to thank Jean-Louis Colliot-Thelene for a very useful remark he made to us at the beginning of this work.

\section{The Fundamental Ideal}
Define $\tilde{W}$ and $\tilde{I}$ to be the sheaves associated to the presheaves $W(-):U\mapsto W(\mathcal{O}_X(U))$ and $I(-):U\mapsto I(W(U))$ respectively.  The sheaves $\tilde{I}^n$ form a filtration of $\tilde{W}$ and the global sections $\mathbb{I}^n:=\Gamma(X,\tilde{I}^n)$ form a filtration of $\Gamma(X,\tilde{W})$.

Assuming that $X$ is proper, smooth, and geometrically connected and that $dim(X)\leq 3$, we first recall a result ~\cite{balmer} that establishes $W(X)=\Gamma(X,\tilde{W})$ and so the $\mathbb{I}^n$ form a filtration of $W(X)$.

\begin{prop}
Let $X$ be a proper, smooth, geometrically connected variety.
$\Gamma(X,\tilde{W})=W(X)$ if dim$X\leq 3$.
\end{prop}

The following result (Proposition \ref{rankbound}) shows that, in order to describe the Witt ring of a variety, it frequently suffices to study the quotients $\mathbb{I}^n/\mathbb{I}^{n+1}$.  Theorem \ref{milnor} then shows that these quotients may be described cohomologically.

\begin{prop}
\label{rankbound}
If $X$ is a proper, smooth, geometrically connected variety with $dim(X)\leq 3$ and the rank of anisotropic forms over $k(X)$ is bounded
then $[E]=0\in W(X)$ if and only if $[E]\in \mathbb{I}^n \: \forall n$.
\end{prop}

\begin{pf}

Let $E\in \mathbb{I}^n$ for all integers $n$.  Then $E_{\eta}\in\tilde{I}^n_{\eta}=I^n(k(X))$.  Hauptsatz X.5.1 of~\cite{lam} shows that the $I^n(k(X))$ are eventually trivial.

\end{pf}

\begin{thm}{(Generalized Theorem of Kerz, Milnor, Orlov, Vishik, Voevodsky et al)}
\label{milnor}

If $X$ is a smooth, quasi-projective variety over an infinite field $k$ then there are isomorphisms
$\xymatrix {\mathscr{K}^n_M/2\mathscr{K}^n_M \ar[r] \ar[d] &  \mathscr{H}^n(X_{et},\mu_2)\\
            \tilde{I}^n/\tilde{I}^{n+1} \ar@{-->}[ur] & }$
\\
\\
where $\mathscr{K}^n_M$ is the Zariski sheafification of Milnor K-theory and
 \noindent $\mathscr{H}^n(\mu_2)$ is the Zariski sheafification of the \'{e}tale cohomology presheaf
\noindent $U\mapsto H^n_{et}(U,\mu_2)$.
\end{thm}

\begin{pf}

These maps are defined on presheaves as follows:\\
$K^n_M/2K^n_M\to H^n(\mu_2):\{a_1,\ldots,a_n\}\mapsto a_1\cup\cdots\cup a_n$\\
$K^n_M/2K^n_M\to I^n/I^{n+1}:\{a_1,\ldots,a_n\}\mapsto<<-a_1,\ldots,-a_n>>$.
That they are isomorphisms at each stalk follows from ~\cite{kerz}.

\end{pf}

We can use this theorem together with the injection $\mathbb{I}^n/\mathbb{I}^{n+1}\hookrightarrow\Gamma(X,\tilde{I}^n/\tilde{I}^{n+1})$ to obtain a description of the first three $\mathbb{I}^n/\mathbb{I}^{n+1}$ (and the images of bilinear spaces therein) in terms of classical invariants.  We first describe the sheaves.

\begin{prop}
Let $X$ be a connected scheme.
Suppose $2$ is a unit on $X$.  Then we have the following.

\begin{enumerate}
	\item $\tilde{W}/\tilde{I}=\mathbb{Z}/2\mathbb{Z}$
	\item $\tilde{I}/\tilde{I}^2=\mathcal{O}^{\times}/\mathcal{O}^{\times 2}$
	\item $\tilde{I}^2/\tilde{I}^3={}_2Br$
\end{enumerate}
\end{prop}

\begin{pf}

	$\tilde{W}/\tilde{I}$ is the sheaf defined on each open $U\subseteq X$ by the presheaf $W(\mathcal{O}_X(U))/I(\mathcal{O}_X(U))$, which is the ring $\mathbb{Z}/2$ for all connected $U$.

	Since the Kummer sequence $0\to\mu_2\to\mathbb{G}_m\to\mathbb{G}_m\to 0$ is exact in the \'{e}tale setting, there is, for each open $U\subseteq X$, a long exact sequence \\
				$\xymatrix{H^0_{et}(U,\mathbb{G}_m)  \ar[r]^{\cdot 2}  & H^0_{et}(U,\mathbb{G}_m) \ar[r]  & H^1_{et}(U,\mu_2) \ar[r]  & H^1_{et}(U,\mathbb{G}_m) \ar[r] & H^1_{et}(U,\mathbb{G}_m)}$	
				
					After the identifications $H^0_{et}(U,\mathbb{G}_m)=\Gamma(U,\mathbb{G}_m)$ and
					$H^1_{et}(U,\mathbb{G}_m)=Pic(U)$, this yields a short exact sequence of presheaves
					
					$0\to\Gamma(U,\mathbb{G}_m)/\Gamma(U,\mathbb{G}_m)^2\to H^1_{et}(U,\mu_2)\to{}_2Pic(U)\to 0$
				
					on each open $U$ so that there is an exact sequence of Zariski sheaves
					
					$0\to\mathcal{O}_X^{\times}/\mathcal{O}_X^{\times 2}\to\mathscr{H}^1_{et}(\mu_2)\to\underline{{}_2Pic(X)}\to 0$.
					
					However, the Picard group is trivial at each stalk, so this becomes a sheaf isomorphism $\mathcal{O}^{\times}/\mathcal{O}^{\times 2}\cong\mathscr{H}^1_{et}(\mu_2)$.
					
	In a similar way, we have a short exact sequence
						
						$0\to Pic(U))/Pic(U))^2\to \mathscr{H}^2_{et}(U,\mu_2)\to{}_2Br(U)\to 0$
						
						and, with no contribution from the Picard group at the level of stalks, we have a sheaf isomorphism $\mathscr{H}^2_{et}(\mu_2){\cong}_2Br$ 

\end{pf}

The following proposition describes the first three quotients $\mathbb{I}^n/\mathbb{I}^{n+1}$ in terms of rank, signed discriminant, and Witt invariant.

\begin{prop}

\begin{enumerate}
			\item The class of $E$ in $W(X)/\mathbb{I}$ is determined by the parity of $n=rank(E)$.
			\item Given an even rank element of $W(X)$, the class of $E\in \mathbb{I}/\mathbb{I}^2$ is determined by the signed discriminant $d_{\pm}(E)=(-1)^{\frac{n(n+1)}{2}}\bigwedge\limits^{rk(E)} E$.
			\item Given an even rank $E\in W(X)$ whose signed discriminant is trivial, the class of E in $\mathbb{I}^2/\mathbb{I}^3$ is determined by the Witt invariant of $E$ (the class of the Clifford algebra of $E$ in the order $2$ Brauer group).
\end{enumerate}

\end{prop}

\begin{pf}

At each stalk,  the composition $\tilde{W}_x\to(\tilde{I}^n/\tilde{I}^{n+1})_x$ ($n\in\{0,1,2\}$) is determined by the invariant specified in the proposition.
\end{pf}

Note that $\mathbb{I}$ consists of the elements of $W(X)$ represented by spaces of even rank and that $\mathbb{I}^2$ consists of those elements of $W(X)$ that are represented by even rank symmetric spaces with trivial signed discriminant.  Furthermore,  every nonzero element of $W(X)/\mathbb{I}$ may be equated in $W(X)/\mathbb{I}$ to any rank $1$ form and every element of $\mathbb{I}/\mathbb{I}^2$ may be equated to a rank $2$ form which shares its signed discriminant and is an orthogonal sum of rank $1$ forms.

\section{Forms of Rank $1$}
\label{rank1}

In order to describe the Witt ring, it will be useful to first address some results and notational issues concerning the part of the Witt ring that is additively generated by forms of rank $1$.  For the purpose of this discussion, we let $X$ be a proper smooth geometrically connected variety over a field k such that $dim(X)\leq 3$ and $W(X)$ be the Witt ring of $X$.  $Q(X)$ will represent the multiplicative subgroup of $W(X)$ represented by rank $1$ forms $(\mathscr{L},\varphi:\mathscr{L}\to\mathscr{L}^{\vee})$ ($(-)^{\vee}$ represents the dual).

\begin{prop}

Given a line bundle $\mathscr{L}\in{ }_2Pic(X)$, the Witt classes $(\mathscr{L},\varphi)$  are in one to one correspondence with the square classes of units of the base field $k$.
\end{prop}

\begin{pf}

First, note that every line bundle $\mathscr{L}$ of order $2$ is equipped with at least one nondegenerate symmmetric bilinear form.

Two forms $(\mathscr{L},\varphi), (\mathscr{L},\psi)\in Q(X)$ may differ at most by a global endomorphism of $\mathscr{L}$ as in the following diagram:

$\xymatrix{\mathscr{L}\ar[r]^{\varphi} \ar[d]_{m_{\ell}=\varphi\psi^{-1}} & \mathscr{L}^{\vee} \\
            \mathscr{L} \ar[ru]_{\psi} & }$
\\
As $\mathscr{E}nd(\mathscr{L})\cong\mathcal{O}_X$, these global endomorphisms are precisely the units of $k$.  Thus, $\varphi$ and $\psi$ differ by multiplication $m_{\ell}$ by some unit $\ell\in k$.

$\varphi$ and $\psi$ represent the same Witt class precisely when $\varphi\cong\psi$ which occurs if and only if $\varphi = m_{\ell}\circ\psi\circ m_{\ell}^{\vee}$ so that  $\varphi=\ell^2\psi$.

This shows that the Witt classes of $X$ associated to $\mathscr{L}$ are in one-to-one correspondence with the square classes of the base field.  Furthermore, the square classes of the base field act multiplicatively to permute the Witt classes associated to $\mathscr{L}$.

\end{pf}

We fix a form on each $\mathscr{L}$ so that the forms on $\mathscr{L}_1,\ldots,\mathscr{L}_n$ multiply to the fixed form on $\mathscr{L}_1\otimes\cdots\otimes\mathscr{L}_n$ and the fixed form on the structure sheaf is the identity.  $<\mathscr{L}>\in W(X)$ will denote the Witt class of a line bundle with such a fixed form and $<s\mathscr{L}>$ will denote multiplication of $<\mathscr{L}>$ by the Witt class $<s>\in W(k)$.  $<s_1\mathscr{L}_1,\ldots,s_n\mathscr{L}_n>$ will denote the orthogonal sum $<s_1\mathscr{L}_1>\perp\cdots\perp<s_n\mathscr{L}_n>$.  The classes $<s>\in Q(k)$ will be identified with the classes $<s\mathcal{O}_X>\in Q(X)$.

We can now give an exact description of $Q(X)$ as a $W(k)$-module.

\begin{prop}

There is a short exact sequence of commutative rings $0\to Q(k)\to Q(X)\to{ }_2Pic(X)\to 0$
\end{prop}

\begin{pf}

Let $<s\mathscr{L}>=<t\mathscr{M}>\in Q(X)$.  We will show that $\mathscr{L}$ and $\mathscr{M}$ are isomorphic as line bundles.

We know that $<s\mathscr{L}>\perp M=<t\mathscr{M}>\perp M'\in Bil(X)$ for metabolic spaces $M$ and $M'$ of the same rank, $2m$.

Taking determinants on both sides, we see that
$$det(<s\mathscr{L}>)det(M)=det(t<\mathscr{M}>)det(M'),$$ $$<s\mathscr{L}>(<-1>^m)=<t\mathscr{M}>(<-1>^m),$$ and $$<s\mathscr{L}>=<t\mathscr{M}>\in Bil(X).$$
\end{pf}

In fact, we have the following.

\begin{prop}

The sub $W(k)$-module $Q(X)$ of $W(X)$ consisting of classes represented by forms on line bundles is isomorphic to\\
$Q(k)\times{ }_2Pic(X)$.
\end{prop}

\section{The Witt Ring of a Curve with Good Reduction over a Non-dyadic Local Field}

If $dim(X)\leq 3$ then the natural map $W(X)\to W(k(X))$ to the Witt ring of the function field is injective ~\cite{balmer}. If $X$ contains a $k$-rational point, the map $W(k)\to W(X)$ induced by the structure map $X\to Spec(k)$ is also injective.  This map gives $W(X)$ the structure of a $W(k)$-module and identifies $W(k)$ with the subring of $W(X)$ generated (as a $W(k)$ module) by classes of forms on $\mathcal{O}_X$.

It is worth noting that the question of whether the composition $W(k)\to W(X)$  is injective is precisely the question of whether the composition $W(k)\to W(k(X))$ is injective, which is in turn equivalent to whether a non-square in the group of units $k^{\times}$ can become a square in $k(X)^{\times}$.

We will now calculate the Witt ring of a smooth geometrically connected projective curve $C$ with good reduction over a non-dyadic local field $K$.  We assume that $C$ contains a $K$-rational point.

We first observe that the elements of the Witt ring $W(C)$ are determined entirely by parity of rank, signed discriminant, and Witt invariant, reducing the calculation of the Witt ring to consideration of classical invariants.

\begin{prop}
\label{higherI}

 $\mathbb{I}^n(C)/\mathbb{I}^{n+1}(C)$ is trivial when $n>2$.

\end{prop}

\begin{pf}

$H^0(C,\mathscr{H}^3(\mu_2))=H^3_{un}(C,\mu_2)=0$ by ~\cite{kato} Proposition 5.2 and the fact that the $2$-cohomological dimension of $C$ is $3$.
\end{pf}

Note that the good reduction assumption is essential here, as examples are known of curves with bad reduction over local fields for which $\mathbb{I}^3/\mathbb{I}^4$ is nontrivial.

The current situation is shown in the following diagram

$\xymatrix{C_k \ar[r] \ar[d] & C_{\nu} \ar[d] & C \ar[l] \ar[d]\\
									Spec(k) \ar[r] & Spec(\mathcal{O}_{\nu}) & Spec(K) \ar[l]}$
\\									
\\
where $C=C_{\nu}\otimes K$ is the generic fiber of the proper, smooth curve $C_{\nu}\to Spec\: \mathcal{O}_{\nu}$ with special fiber $C_k=C_{\nu}\otimes k$ smooth over $Spec(k)$. The $K$ point of $C$ determines an $\mathcal{O}_{\nu}$ rational point of $C_{\nu}$ by the valuative criterion of properness.

A key tool is the cohomological calculation  of these invariants and so we need to understand the relationship between them for the curve $C_{\nu}$ and $C_k$.

\begin{lem}

\label{basechange}

Let $f:X\rightarrow Spec\left( \mathcal{O}_{\nu }\right) $ be a proper
morphism over the Hensel ring $\mathcal{O}_{\nu }$ with residue field $k$
and let $X_{k}$ denote the fibre of $p$ over $Spec\left( k\right) .$ Then
the base change map $H^{p}\left( X,\mu _{n}\right) \rightarrow H^{p}\left(
X_{k},\mu _{n}\right) $ is an isomorphism.
\end{lem}

\begin{pf}
If $\mathcal{O}_{\nu }$ is strictly local, this is just the proper base
change isomorphism. Since finite, \'{e}tale covers of $Spec\left( k\right) $
and $Spec\left( \mathcal{O}_{\nu }\right) $ are the same, the result follows
immediately from the Hochschild-Serre spectral sequences.
\end{pf}

Once again, the forms on an individual line bundle are in one to one correspondence with the rank 1 forms in the Witt ring of the base field, of which there are four.  We will equip each line bundle with a fixed form as described in Section ~\ref{rank1} which extends to a non-degenerate form on $C_k$.

\begin{prop}
${}_2Pic(C_k)={}_2Pic (C_{\nu})={}_2Pic (C)$
\end{prop}

\begin{pf}

Since $C_{\nu }$ is proper and smooth over $\mathcal{O}_{\nu }$ with a
rational point, its Picard scheme has an abelian variety $Pic_{C_{\nu }/%
\mathcal{O}\nu }^{0}$ as its connected component and, for any $\mathcal{O}%
_{\nu }$ scheme $T,$ $Pic_{C_{\nu }/\mathcal{O}\nu }^{0}\left( T\right)
=Pic^{0}\left( C\times _{\mathcal{O}_{\nu }}T\right) /Pic\left( T\right) $
where  $Pic^{0}\left( C\times _{\mathcal{O}_{\nu }}T\right) $ consists of
line bundles on $C\times _{\mathcal{O}_{\nu }}T$ of degree $0$ on each fibre
over $T.$ In particular \thinspace $_{2}Pic^{0}\left( C_{k}\right)
=\thinspace _{2}Pic_{C_{\nu }/\mathcal{O}\nu }^{0}\left( k\right) ,$ \thinspace $%
_{2}Pic^{0}\left( C_{\nu }\right) =\thinspace _{2}Pic_{C_{\nu }/\mathcal{O}\nu }^{0}\left(
\mathcal{O}_{\nu }\right) $ and \thinspace $_{2}Pic^{0}\left( C\right)
=\thinspace _{2}Pic_{C_{\nu }/\mathcal{O}\nu }^{0}\left( K\right) .$ But $\,_{2}Pic_{C_{\nu
}/\mathcal{O}\nu }^{0}$ is a finite \'etale group scheme since $K$ is
non-dyadic and $\mathcal{O}_{\nu }$ is henselian. Consequently $%
_{2}Pic^{0}\left( C_{k}\right) =\,_{2}Pic^{0}\left( C_{\nu }\right) .$
Moreover the valuative criterion for separation and properness shows that $%
\,_{2}Pic^{0}\left( C_{\nu }\right) =\,_{2}Pic^{0}\left( C\right) .$

\end{pf}

Thus, associated to each line bundle $\mathscr{L}\in{}_2Pic(C)$, there is a unique line bundle $\mathscr{L}_{\nu}\in{}_2Pic(C_{\nu})$.  We select a Witt class represented by $(\mathscr{L}_{\nu},\varphi)$ and note that $(\mathscr{L},\varphi\otimes 1)$ represents a Witt class of $C$.  Furthermore, $(\mathscr{\bar{L}},\bar{\varphi})$ represents a Witt class of $C_k$ and every rank $1$ element of $W(C_k)$ may be obtained by base change from a form over $C_{\nu}$.  Note that the selected forms multiply as expected and we may take the form on $\mathcal{O}_C$ to be $<1>=<1\mathcal{O}_C>$.  There are four Witt classes associated to each $\mathscr{L}$.  These are $<\mathscr{L}>$, $<s\mathscr{L}>$, $<\pi\mathscr{L}>$, and $<s\pi\mathscr{L}>$ where $<1>\neq<s>\in W(k)$ and $\pi$ is the uniformizing parameter of $\mathcal{O}_{\nu}$.

Since the forms $<\mathscr{L}>$ and $<s\mathscr{L}>$ arise from forms over $C_k$, these forms behave like forms over a curve over a finite field.  In particular, we have relations $$<u\mathscr{L},v\mathscr{M}>=<1,uv\mathscr{L}\mathscr{M}>$$ and
$$<\pi u\mathscr{L},\pi v\mathscr{M}>=<\pi,\pi uv\mathscr{L}\mathscr{M}>$$ in $W(C)$ where $u,v\in k^{\times}/k^{\times 2}$~\cite{funk}.

We can equate each form to a rank $1$ form in $W(C)/\mathbb{I}$ and each even dimensional form to a sum of rank $1$ forms in $\mathbb{I}/\mathbb{I}^2$.  It remains to describe $\mathbb{I}^2/\mathbb{I}^3$, which we do by describing the order $2$ part of the Brauer group.

Given a square class, $s$, of $K^{\times}$ and an order $2$ line bundle, $\mathscr{L}$, we will construct an Azumaya algebra $\big(\frac{s\mathscr{L},\pi}{C}\big)$ modeled after the quaternion algebras.
Let $\mathscr{L}\in{ }_2Pic(C)$, $<s>\in W(k)$.
Take $f:\mathscr{L}\otimes\mathscr{L}\to\mathcal{O}_C$ to be the isomorphism corresponding to the form on $<s>\mathscr{L}$.  Let multiplication of graded rings on $\mathcal{O}_C\oplus\mathscr{L}$ be induced by $f$ (locally $(a,\ell)\cdot(a',\ell')=(aa'+f(\ell,\ell'),a\ell'+a'\ell)$) and on $\mathcal{O}_C\oplus\mathcal{O}$ by the isomorphism $\mathcal{O}\otimes\mathcal{O}\to\mathcal{O}$ induced by the form $<\pi>$ (locally $(a,b)\cdot(a',b')=(aa'+\pi^2 bb',ab'+a'b)$).

Define $\big(\frac{s\mathscr{L},\pi}{C}\big)$ to be the graded tensor product $(\mathcal{O}_C\oplus\mathscr{L})\hat{\otimes}(\mathcal{O}_C\oplus\mathcal{O})$

Note that,  as an $\mathcal{O}_C$-module, $\big(\frac{s\mathscr{L},\pi}{C}\big)=\mathcal{O}_C\oplus\mathscr{L}\oplus\mathscr{L}\oplus\mathcal{O}_C$ and that multiplication on the second copy of $\mathscr{L}$ is given by $\pi f$.  Locally this is just the quaternion algebra $\big(\frac{s\xi,\pi}{\mathcal{O}_x}\big)$ additively generated by $i, j, k$ with $i^2=\xi$, $j^2=\pi$, and $ij=-ji=k$ where $<\xi>$ is the image of $\mathscr{L}$ in $W(\mathcal{O}_x)$. In particular $\big(\frac{s\mathscr{L},\pi}{C}\big)$ is a sheaf of Azumaya algebras on $C$. The norm form of $\big(\frac{\xi,\pi}{k(C)}\big)$ is $<1,-s\xi,-\pi,s\xi\pi>$, which is the image in $W(k(C))$ of the form $<1,-s\mathscr{L},-\pi,s\pi\mathscr{L}>$.

These $\big(\frac{s\mathscr{L},\pi}{C}\big)$ will comprise the order $2$ Brauer group, as shown in the following two theorems.



\begin{thm}
Let $C$ be a curve with good reduction over a non-dyadic local field $K$
that has a $K$ rational point $P$.
The order 2 Brauer group, ${}_2Br(C)$, fits into a short exact sequence
$0\to k^{\times}/k^{\times 2}\to{}_2Br(C)\to{}_2Pic(C)\to 0$.
\end{thm}

\begin{pf}
Gabber's proof of purity \cite{Gabber} extends Corollary 6.2 \cite{GBIII} to
our setting. The localization sequence, suitably interpreted, thus becomes
the exact sequence
\[
0\rightarrow Br\left( C_{\nu }\right) \left( \ell \right) \rightarrow
Br\left( C\right) \left( \ell \right) \rightarrow H^{1}\left( C_{k},\mathbb{Q%
}_{\ell }/\mathbb{Z}_{\ell }\right) \rightarrow H^{3}\left( C_{\nu
},G_{m}\right) \left( \ell \right) \overset{i}{\rightarrow }H^{3}\left(
C,G_{m}\right) \left( \ell \right)
\]%
for any prime $\ell \neq 2$. By Artin's result \cite[Theorem 3.1]{GBIII} and
the vanishing of $Br\left( C_{k}\right) $, $Br\left( C_{\nu }\right) =0.$ We
will show that the map $i$ is
injective for $\ell \neq p=char\left( k\right) .$ Then $Br\left( C\right)
\left( 2 \right) \approxeq H^{1}\left( C_{k},\mathbb{Q}_{2}/\mathbb{Z}%
_{ 2 }\right) $ and so \thinspace $_{2}Br\left( C\right) \approxeq
H^{1}\left( C_{k},\mathbb{Z}/2\right) $ as desired.

First we need to calculate $H^{3}\left( C_{\nu },G_{m}\right) \left( \ell
\right) $. So take the limit over the Kummer sequences for powers of $\ell $
and observe that $H^{2}\left( C_{\nu },G_{m}\right) =0.$
Then we find that $H^{3}\left( C_{\nu },\mu _{\ell ^{\infty }}\right)
\approxeq H^{3}\left( C_{\nu },G_{m}\right) \left( \ell \right) .$ Similarly
$H^{3}\left( C_{k},\mu _{\ell ^{\infty }}\right) \approxeq H^{3}\left(
C_{k},G_{m}\right) \left( \ell \right) .$
Now Lemma \ref{basechange}
establishes the base
change isomorphism $H^{3}\left( C_{\nu },\mu _{\ell ^{\infty }}\right)
\approxeq H^{3}\left( C_{k},\mu _{\ell ^{\infty }}\right) $. But the
Hochschild-Serre spectral sequence $H^{p}\left( G\left( k_{s}/k\right)
,H^{q}\left( C_{k,s},G_{m}\right) \right) \Longrightarrow H^{n}\left(
C_{k},G_{m}\right) $ and Tsen's theorem establish the isomorphism $%
H^{2}\left( G\left( k_{s}/k\right) ,H^{1}\left( C_{k,s},G_{m}\right) \right)
\approxeq H^{3}\left( C_{k},G_{m}\right)$ where $C_{k,s}$ is the base extension
of $C_{k}$ to the separable closure of $k$. But $C_{k}$ has a $k$ rational
point and so \ there is a $G\left( k_{s}/k\right) $ isomorphism $H^{1}\left(
C_{k,s},G_{m}\right) =Pic^{0}\left( C_{k,s}\right) \oplus \mathbb{Z}$ where $%
\mathbb{Z}$ has trivial $G\left( k_{s}/k\right) $ action. $Pic^{0}\left(
C_{k,s}\right) $ is torsion and so has trivial $G\left( k_{s}/k\right) $
second cohomology while $\mathbb{Z}$ fits into a short exact sequence of
trivial $G\left( k_{s}/k\right) $ modules, $0\rightarrow \mathbb{Z}%
\rightarrow \mathbb{Q}\rightarrow \mathbb{Q}/\mathbb{Z}\rightarrow 0.$ So $%
H^{3}\left( C_{k},G_{m}\right) =Hom\left( G\left( k_{s}/k\right) ,\mathbb{Q}/%
\mathbb{Z}\right) \approxeq \mathbb{Q}/\mathbb{Z}.$ (For $\ell \neq
p=char\left( k\right)$, we could also use the Kummer sequence for $\ell ^{n}
$ to identify $H^{3}\left( C_{k},\mu _{\ell ^{n}}\right) =\mathbb{Z}/\ell
^{n}$ with $\,_{\ell ^{n}}H^{3}\left( C_{k,}G_{m}\right) $ and then observe
that  $\mathbb{Z}/\ell ^{n}=\,_{\ell ^{n}}\mathbb{Q}/\mathbb{Z}$ maps into
the summand of $H^{2}\left( G\left( k_{s}/k\right) ,H^{1}\left(
C_{k,s},G_{m}\right) \right) $ defined by the summand $\mathbb{Z}%
\hookrightarrow Pic\left( C_{k,s}\right) $ since the identification comes
from $H^{2}\left( C_{k,s},\mu _{\ell ^{n}}\right) =\mathbb{Z}/\ell ^{n}\left[
P_{k}\right] $ where $P_{k}$,  the $k$ rational point on $C_{k}$, is the
reduction of the $\mathcal{O}_{\nu }$ point $P_{\nu }$ which, in turn is the
extension of the $K$ rational point $P$ on $C$ to an $\mathcal{O}_{\nu }$
point on $C_{\nu }$.)

Finally we show that $i$ is injective by noting that the existence of $P$ on
$C$ that reduces to $P_{k}$ on $C_{k}$ extends to an identification of the
generator of $\mathbb{Z}=Pic\left( C\right) /Pic^{0}\left( C\right) $ with a
generator of $Pic\left( C_{k}\right) /Pic^{0}\left( C_{k}\right) .$
Consequently
\[
H^{2}\left( G\left( k_{s}/k\right) ,H^{1}\left( C_{k,s},G_{m}\right) \right)
=H^{2}\left( G\left( k_{s}/k\right) ,H^{1}\left( C_{\nu }^{s},G_{m}\right)
\right) \rightarrow H^{2}\left( G\left( K_{s}/K\right) ,H^{1}\left(
C_{s},G_{m}\right) \right)
\]%
may be interpreted as $Hom\left( G\left( k_{s}/k\right) ,\mathbb{Q}/\mathbb{Z%
}\right) \hookrightarrow Hom\left( G\left( K_{s}/K\right) ,\mathbb{Q}/%
\mathbb{Z}\right)$. (Here $C_{\nu }^{s}$ is the extension of $C_{\nu}$ to the strict henselization of $\mathcal{O}_{\nu }$.) But, in the Hochschild-Serre spectral sequence over $K,$
$H^{p}\left( G\left( K_{s}/K\right) ,H^{q}\left( C_{s},G_{m}\right) \right)
\Longrightarrow H^{n}\left( C,G_{m}\right) ,$ $H^{2}\left( G\left(
K_{s}/K\right) ,\text{ }H^{1}\left( C_{s},G_{m}\right) \right) \subseteq
H^{3}\left( C,G_{m}\right) $ because $E_{\infty }^{2,1}=E_{2}^{2,1}$ as we observe by noting $%
cd\left( K\right) =2$ and $Br\left( C_{s}\right) =0.$
\end{pf}
The final step is to explicitly identify the algebras that lie in ${}_2Br(C)$
\begin{thm}

The order $2$ Brauer group of a smooth curve with good reduction over a local field consists of the distinct quaternions $\big(\frac{s\mathscr{L},\pi}{C}\big)$ where $s\in k^{\times}/k^{\times 2}$ and $\mathscr{L}\in{}_2Pic(C)$.
\end{thm}
\begin{pf}

It suffices to show that the quaternions are distinct.  That they comprise the entire order $2$ Brauer group follows by a cardinality argument.   We work over the function field via the injections ${}_2Br(C)\hookrightarrow{}_2Br(k(C))$ and $W(C)\hookrightarrow W(k(C))$, taking $<s\xi>$ and $<t\zeta>$ to be the images of $<s\mathscr{L}>$ and $<t\mathscr{M}>$ in $W(k(C))$.

If $\big(\frac{s\xi,\pi}{k(C)}\big) = \big(\frac{t\zeta,\pi}{k(C)}\big)$ then their norm forms are equal in $W(k(C))$ and we have a relation

$<1,-\pi,-\xi,\pi\xi,-1,\pi,\zeta,-\pi\zeta>=<1,-\pi,-\xi\zeta,\pi\xi\zeta>=0\in W(k(C))$.

This is the norm form of $\big(\frac{st\xi\zeta,\pi}{k(C)}\big)$.

Thus, it suffices to show that $\big(\frac{s\xi,\pi}{k(C)}\big)$ is nontrivial when $<\xi>\neq<1>$.

If $\big(\frac{s\xi,\pi}{k(C)}\big)$ is trivial in ${}_2(Br(k(C))$ then there is an identity $a'^2-b'^2\xi=\pi$ in $k(C)$ (~\cite{lam} Theorem III.2.7).

If we take $p$ to be the extension of the $K$-rational point over the closed point of $Spec(\mathcal{O}_{\nu})$ then the following are true of the local ring $\mathcal{O}_{C_{\nu},p}$:

\begin{itemize}
	\item $\mathcal{O}_{C_{\nu},p}$ is a unique factorization domain
	\item $\pi$ is not a unit of $\mathcal{O}_{C_{\nu},p}$
	\item $\pi$ is a prime/irreducible element of $\mathcal{O}_{C_{\nu},p}$
	\item The field of fractions of $\mathcal{O}_{C_{\nu},p}$ is $k(C)$
\end{itemize}

This means that we can write $a'=\frac{a}{d}$ and $b'=\frac{b}{e}$ where $a,b,d,e\in\mathcal{O}_{C_{\nu},p}$ and $(a,d)=(b,e)=1$, giving a relation $a^2e^2-b^2d^2\xi=\pi d^2e^2$ in $\mathcal{O}_{C_{\nu},p}$.  We can then factor all copies of $\pi$ from both sides of this equation and reduce modulo $\pi$ to get a relation $\bar{a}^2-\bar{b}^2s\xi=0$ in the local ring $\mathcal{O}_{C_k,p}$.  This means that $<s\xi>=<1>$ at the stalk $\mathcal{O}_{C_k,p}$ and in the function field of $C_k$.


\end{pf}

We can now calculate the Witt ring of $C$ by describing its unique representatives.

\begin{thm}
\label{thm:curvelocalfield}

Each nontrivial element of $W(C)$ has a representative of one of the following forms:
$<s\mathscr{L}>$, $<t \pi\mathscr{M} >$, $<1,s\mathscr{L}>$, $<s\mathscr{L},t \pi\mathscr{M}>$, $<\pi,t \pi\mathscr{M}>$, $<1,s\mathscr{L},t \pi\mathscr{M}>$, $<s\mathscr{L},\pi,t\pi\mathscr{M}>$, or $<1,s\mathscr{L},\pi,t \pi\mathscr{M}>$ where $<s>,<t>\in W(k)$ and $\mathscr{L},\mathscr{M}\in{ }_2Pic(C)$.

Furthermore, these forms represent distinct Witt classes provided the fewest possible line bundles are used to represent each class and even rank forms are nontrivial except in the obvious cases.
\end{thm}

\begin{pf}

Let $n$ be the cardinality $|{ }_2Pic(C)|$ 

$W(C)/\mathbb{I}$ is the two element group with nontrivial representative $<1>$, \\
$\mathbb{I}/\mathbb{I}^2$ has $4n$ elements represented by forms $<-1,\delta>$, and $\mathbb{I}^2/\mathbb{I}^3$ has $2n$ elements represented by quaternion norm forms $<1,-s\mathscr{L},-\pi,s\pi\mathscr{L}>$.

This gives us two important facts.

$W(C)$ is a ring with $16n^2$ elements.




Every nontrivial element of $W(C)$ may be written\\
 $<1>^{d_1}\perp<-1,d>^{d_2}\perp<1,-s\mathscr{L},-\pi,s\pi\mathscr{L}>^{d_3}$, $d_i\in\{0,1\}$ (using Proposition~\ref{higherI}).

We will now perform the calculation to show that \\
$<1>\perp<-1,\delta>\perp<1,-s\mathscr{L},-\pi,s\pi\mathscr{L}>$ has a representative of one of the listed forms.

$E=<1>\perp<-1,\delta>\perp<1,-s\mathscr{L},-\pi s\mathscr{L},\pi>=\\
<1,-1,\delta,1,-s\mathscr{L},-\pi s\mathscr{L},\pi>=<1,-s\mathscr{L},\delta,\pi,-\pi s\mathscr{L}>$.

If $\delta=<t\mathscr{M}>$ is a form which extends to a non-degenerate form over $C_{\nu}$ then $<1,-s\mathscr{L},\delta,\pi,-\pi s\mathscr{L}>=<st\mathscr{L}\mathscr{M},\pi,-s \pi\mathscr{L}>$

If $\delta=<t\pi\mathscr{M}>$ is a form which does not extend to a non-degenerate form over $C_{\nu}$ then $<1,-s\mathscr{L},\delta,\pi,-\pi s\mathscr{L}>=<1,-s\mathscr{L},st\pi\mathscr{L}\mathscr{M}>$.

This completes the case $d_1=d_2=d_3=1$.  The calculations for other values of $d_i$ are analagous or obvious.

It remains to show that,assuming the fewest number of line bundles possible are used to represent each form, the forms described are distinct.

We know that $W(C)$ has $16n^2-1$ nontrivial elements and count the forms that cannot be written in an obvious way as a sum of fewer line bundles.  There are $2n$ forms each $<s\mathscr{L}>$ and $<t\pi\mathscr{M}>$, $2n-1$ forms each\\
 $<1,s\mathscr{L}>$ and $<\pi,t\pi\mathscr{M}>$, $(2n)^2$ forms $<s\mathscr{L},t\pi\mathscr{M}>$, $2n(2n-1)$ forms each $<1,s\mathscr{L},t\pi\mathscr{M}>$ and $<s\mathscr{L},\pi,t\pi\mathscr{M}>$, and $(2n)^2-2n-2n+1$ forms $<1,s\mathscr{L},\pi,t\pi\mathscr{M}>$.  This is a total of $16n^2-1$ forms.

A cardinality argument shows that the forms of this proposition are distinct.

\end{pf}

We conclude this paper with two theorems describing the algebraic structure of $W(C)$ as it relates to $W(C_k)$.  This generalizes the analogous result for the Witt ring of a local field (~\cite{lam} Section VI.1).

\begin{thm}
\label{thm:curvelocal2}

$W(C)$ is an abelian group isomorphic to $W(C_k)\oplus W(C_k)$.
\end{thm}

\begin{pf}

There is an inclusion $W(C_k)\to W(C):\bar{\mathscr{L}} \mapsto\mathscr{L}$ which maps $W(C_k)$ onto the elements of the forms $<s\mathscr{L}>$ or $<1,s\mathscr{L}>$ in $W(C)$.  We identify $W(C_k)$ with its image under this inclusion.  The cokernel of this map is $<\pi>W(C_k)$ and $<\pi>W(C_k)\cong W(C_k)$ by factoring out $<\pi>$.  This gives a short exact sequence.

$0\to W(C_k)\to W(C)\to W(C_k)\to 0$

We will define a left splitting map $W(C)\to W(C_k)$.

$<1,-\pi>W(C_k)$ is an ideal in $W(C)$ as\\
 $<s\mathscr{L}><1,-\pi>\in<1,-\pi>W(C_k)$ and\\
 $<s\pi\mathscr{L}><1,-\pi>=<s\mathscr{L}><\pi,-1>=<-s\mathscr{L}><1,-\pi>$, which is also in $<1,\pi>W(C_k)$.  Furthermore, the composition $f$ in the following diagram is an isomorphism.

$\xymatrix{W(C_k) \ar[r]^i \ar[rd]_f& W(C) \ar[d]^p \\
									& W(C)/<1,-\pi>W(C_k)}$
\\
\\
is an isomorphism.

The desired splitting map is $f^{-1}\circ p$.
\end{pf}

Note that the splitting map $f^{-1}\circ p$ is a map of commutative rings.

\begin{thm}

$W(C)$ is isomorphic as a ring to the group ring $W(C_k)[G]$ where $G$ is a two element group.
\end{thm}

\begin{pf}

Taking $<1>$ and $<\pi>$ as the representatives for $G$, this is a corollary of theorem~\ref{thm:curvelocal2}.
\end{pf}

\bibliographystyle{plain}
\bibliography{bib1}

\begin{thebibliography}{1}

\bibitem{auel}
A.~Auel.
\newblock Remarks on the milnor conjecture over schemes.
\newblock {\em Arxiv preprint arXiv:1109.3294}, 2011.

\bibitem{balmer}
P.~Balmer and C.~Walter.
\newblock A gersten-witt spectral sequence for regular schemes.
\newblock {\em Annales Scientifiques de l'\'Ecole Normale Sup\'erieure},
  35(1):127--152, 2002.

\bibitem{Gabber}
Kazuhiro Fujiwara.
\newblock A proof of the absolute purity conjecture (after {G}abber).
\newblock In {\em Algebraic geometry 2000, {A}zumino ({H}otaka)}, volume~36 of
  {\em Adv. Stud. Pure Math.}, pages 153--183. Math. Soc. Japan, Tokyo, 2002.

\bibitem{funk}
J.~Funk and R.~Hoobler.
\newblock The witt ring of a smooth projective curve over a finite field.
\newblock {\em Arxiv preprint arXiv:1210.3109}, 2012.

\bibitem{GBIII}
Alexander Grothendieck.
\newblock Le groupe de {B}rauer. {III}. {E}xemples et compl\'ements.
\newblock In {\em Dix {E}xpos\'es sur la {C}ohomologie des {S}ch\'emas}, pages
  88--188. North-Holland, Amsterdam, 1968.

\bibitem{kato}
Kazuya Kato.
\newblock A {H}asse principle for two-dimensional global fields.
\newblock {\em J. Reine Angew. Math.}, 366:142--183, 1986.
\newblock With an appendix by Jean-Louis Colliot-Th{\'e}l{\`e}ne.

\bibitem{kerz}
Moritz Kerz.
\newblock The {G}ersten conjecture for {M}ilnor {$K$}-theory.
\newblock {\em Invent. Math.}, 175(1):1--33, 2009.

\bibitem{lam}
T.~Y. Lam.
\newblock {\em Introduction to quadratic forms over fields}, volume~67 of {\em
  Graduate Studies in Mathematics}.
\newblock American Mathematical Society, Providence, RI, 2005.

\end{thebibliography}

\end{document}